\documentclass[11pt, oneside]{article}%{amsart} 	% use "amsart" instead of "article" for AMSLaTeX format
\usepackage{geometry}                		% See geometry.pdf to learn the layout options. There are lots.
\usepackage{graphicx}				% Use pdf, png, jpg, or eps� with pdflatex; use eps in DVI mode
\usepackage{subcaption} % For subfigures
\usepackage{amssymb}
\usepackage{float}
\usepackage{xcolor}% colored fonts

\usepackage[english]{babel}
\usepackage{amsthm}    %per i teoremi
\usepackage{amsfonts}  %per il carattere
\usepackage{mathrsfs}  %per il carattere
\usepackage{amsmath}   %per scrivere in grassetto simboli matematici
\usepackage{verbatim}

\usepackage{amsmath,amsfonts,amssymb,amsthm,epsfig,epstopdf,titling,url,array}

\theoremstyle{plain}

\theoremstyle{definition}
\newtheorem{defn}{Definition}[section]

\theoremstyle{remark}
\newtheorem{rem}{Remark}

\begin{document}

\title{ Deformed Boson Algebras and $\mathcal{W}_{\alpha,\beta,\nu}$-Coherent States: A New Quantum Framework}
\author{R. Droghei$^1$\\
CNR-ISMAR, Institute of Marine Sciences, 00133 Rome, Italy\\
E-Mail Address: riccardo.droghei@artov.ismar.cnr.it }

\maketitle

\begin{abstract}
We introduce a novel class of coherent states, termed $\mathcal{W}^{(\bar{\alpha},\bar{\nu})}(z)$-coherent states, constructed using a deformed boson algebra based on the generalised factorial $[n]_{\alpha,\beta,\nu}!$. This algebra extends conventional factorials, incorporating advanced special functions such as the Mittag-Leffler and Wright functions, enabling the exploration of a broader class of quantum states. The mathematical properties of these states, including their continuity, completeness, and quantum fluctuations, are analysed. A key aspect of this work is the resolution of the Stieltjes moment problem associated with these states, achieved through the inverse Mellin transformation method. The framework provides insights into the interplay between the classical and quantum regimes, with potential applications in quantum optics and fractional quantum mechanics. By extending the theoretical landscape of coherent states, this study opens avenues for further exploration in mathematical physics and quantum technologies. 

\end{abstract}
\smallskip 


%---------------------------------------------------------------------------------------------------------------------------------------------------------------------------------------

\section{Introduction}\label{sec:1}

\setcounter{section}{1} \setcounter{equation}{0}

The coherent states (CS) were introduced by Erwin Schr\"{o}dinger in the far 1926 \cite{S26} in the study of the harmonic oscillator. Still, they only recently were in use. In fact, Glauber (\cite{Glaub63}) was the first to use the name "coherent states" in the field of quantum optics. CS are quantum states that provide a strict relationship between classical and quantum behaviour.
Since the early days of quantum mechanics, various deformations and generalisations of CS and canonical commutation relations have been proposed. The construction of generalised CS through the solution of Stieltjes moment problems has been extensively discussed in the literature \cite{Klauder63}, \cite{KPS2001}, \cite{GBG2010}, \cite{KS91}.

{A significant motivation for these extensions arises from the incorporation of fractional calculus, which allows for a more generalized and flexible approach to quantum mechanics. In recent years, fractional differential operators have been used to model non-local and memory-dependent effects, which naturally emerge in complex quantum systems, such as anomalous diffusion and generalized uncertainty principles \cite{Herrmann2011}. Inspired by these developments, we extend the traditional boson algebra framework by introducing a new class of coherent states, the $\mathcal{W}^{(\bar{\alpha},\bar{\nu})}(z)$-coherent states, derived from a deformed boson algebra incorporating special functions from fractional calculus, such as Mittag-Leffler and Wright functions. These generalizations provide a richer structure for quantum states, offering deeper insights into quantum fluctuations, uncertainty relations, and their connection to fractional quantum mechanics. }

An ensemble of states, in the Dirac notation $|z\rangle$, where $z$ is an element of an appropriate space endowed with the notion of continuity, is called a set of CS if it has the following two properties.

The first property is the continuity: the vector $|z>$ is a strongly continuous function of the label $z$, i.e.
$$|z-z'|\rightarrow0\,\implies \||z\rangle-|z'\rangle\|\rightarrow0$$ where $\||\psi\rangle\|\equiv\langle\psi|\psi\rangle^{1/2}$.

The second property is the completeness (resolution of unity); there exists a positive function $W(|z|^2)$ such that the unity operator $\hat{I}$ admits the "resolution of unity":
\begin{equation}\label{completeness}
\iint_\mathbb{C} d^2z|z\rangle W(|z|^2)\langle z|=\hat{I}=\sum_{n=0}^\infty|n\rangle\langle n|
\end{equation}
where $|n\rangle$ is a {complete} set of orthonormal eigenfunctions of a Hermitian operator.
However, while continuity in $z$ is easy to verify, the condition in Eq. (\ref{completeness}) imposes a significant restriction on the choice of parameters in the definition of CS {(we will see that for the deformed CS introduced in the paper, the parameters satisfy { $0<\alpha+\beta\le 2$} and $\nu>\alpha-1$)}. Only a relatively small number of distinct sets of CS are known for which the function $ W(|z|^2)$ can be explicitly determined. As a result, the family of CS remains limited in size \cite{KS85}.

During the past thirty years, progress has been made in resolving the unity condition for selected parameter choices in CS \cite{PS99}, \cite{SP20002001}, \cite{SPS99}, \cite{PBDHS2010}.
%-------------------------------------------------------------------------------------------------------------------------------------------------------------
{Understanding the physical significance of deformation parameters is crucial for further development. Various generalizations of boson algebras have found applications in quantum superintegrable systems (see \cite{BDK94}), nuclear physics with q-deformed bosons \cite{DGI2001}, and algebraic methods for q-deformed many-body systems \cite{GGR95}. These contributions highlight the significance of deformed boson algebras in different physical contexts. Moreover, previous studies have explored f-deformed coherent states and their applications \cite{TP2007}, \cite{RT2004}, \cite{MMSZ97}. Deformation schemes have also been investigated in the algebras SU(1,1) and SU(2) \cite{SIJ2009}, \cite{BDK93}, further demonstrating the breadth of research in this field. In light of these studies, our work aims to extend the existing framework by introducing a new class of deformed coherent states, emphasising their mathematical and physical implications.}
%-------------------------------------------------------------------------------------------------------------------------------------------------------------

The physical motivation behind the structure of CS is to propose a general linear combination of basis states $|n\rangle$, with coefficients specifically designed to satisfy the equation (\ref{completeness}). These coefficients can often be linked to a specific Hamiltonian $\hat{H}\neq\hat{H}_0$, where $\hat{H}_0$ is the Hamiltonian of the linear harmonic oscillator. As we demonstrate below, a relatively general class of CS is associated with the special function of three parameters $\mathcal{W}_{\alpha,\beta,\nu}$, for which the above conditions can be satisfied. In particular, for different parameter sets $(\alpha, \beta,\nu)$, the explicit form of $ W(|z|^2)$ is derived.

The structure of this paper is as follows: in Sec. 2, we introduce the deformed boson algebra and define the $\mathcal{W}_{\alpha,\beta,\nu}$-coherent states, discussing their mathematical properties, including continuity and completeness. Sec. 3 examines their physical implications, particularly quantum fluctuations and their connection to generalised uncertainty relations. In Sect. 4, we address the Mandel parameter to explore the nonclassical nature of these states. Finally, in the Appendix, we provide detailed asymptotic analyses and verify the conditions under which the resolution of unity is satisfied.

%---------------------------------------------------------------------------------------------------------------------------------------------------------------------------------------

%---------------------------------------------------------------------------------------------------------------------------------------------------------------------------------------

\section{Deformed Boson Algebra and $\mathcal{W}_{\alpha,\beta,\nu}$-coherent states}\label{sec:2}

\setcounter{section}{2} \setcounter{equation}{0}

This section introduces a deformed boson algebra and a novel set of generalised CS, presenting a deformed factorial, denoted $[n]_{\alpha,\beta,\nu}!$, {(where $(\alpha\in [0,1],\beta\in (0,1])$ and $\nu > \alpha-1$)}, which generalizes the classical factorial function. The generalised factorial is defined as follows:

\begin{equation}\label{genfact}
[n]_{\alpha,\beta,\nu}!=[1]_{\alpha,\beta,\nu}...[n]_{\alpha,\beta,\nu}=\prod_{i=1}^n\frac{\Gamma(\beta i+1)}{\Gamma(\beta i+1-\alpha)}\frac{\Gamma(\beta n+1-\alpha+\nu)}{\Gamma(1-\alpha+\nu)}
\end{equation}

 where the box function $[n]_{\alpha,\beta,\nu}$, which determines the deformation, is defined as follows.

\begin{equation}
[n]_{\alpha,\beta,\nu}=\frac{\Gamma(\beta n+1)}{\Gamma(\beta n+1-\alpha)}\frac{\Gamma(\beta n+1-\alpha+\nu)}{\Gamma(\beta (n-1)+1-\alpha+\nu)};
\end{equation} 
{if $n\in\mathbb{N}$, and $[0]_{\alpha,\beta,\nu}=0$.}
 It is simple to verify, also using the telescoping product, that $[0]_{\alpha,\beta,\nu}!=1$ and $[n]_{\alpha,\beta,\nu}!=[n]_{\alpha,\beta,\nu}\cdot[n-1]_{\alpha,\beta,\nu}!$. 
 This generalised factorial allows for a more flexible and extended form of the factorial function, suitable for various applications in advanced mathematical contexts such as combinatorics, special functions, and complex analysis. 

\begin{rem}
Setting the parameter of the generalised factorial (\ref{genfact}), we can obtain the following particular cases attributable to well-known special functions.

{The simplest case $[n]_{0,1,0}!=n!$ is reducible to the classical boson algebra. Furthermore, $[n]_{0,\alpha,\beta-1}!=\frac{\Gamma(\alpha n+\beta)}{\Gamma(\beta)}$ is related to the Mittag-Leffler function investigated in \cite{SPS99}, and $[n]_{1,\lambda,\mu}!=\frac{\lambda^n}{\Gamma(\mu)}n!\Gamma(\lambda n+\mu)$ is related to the Wright function investigated in \cite{GGM19} and \cite{GM23}.}
\end{rem}

Now, we define the $\mathcal{W}_{\alpha,\beta,\nu}$-deformed boson algebra generated by the set of operators. 

\begin{figure}[ht]
    \centering
    % First image
    \begin{subfigure}[b]{0.45\linewidth}
        \centering
        \includegraphics[width=\linewidth]{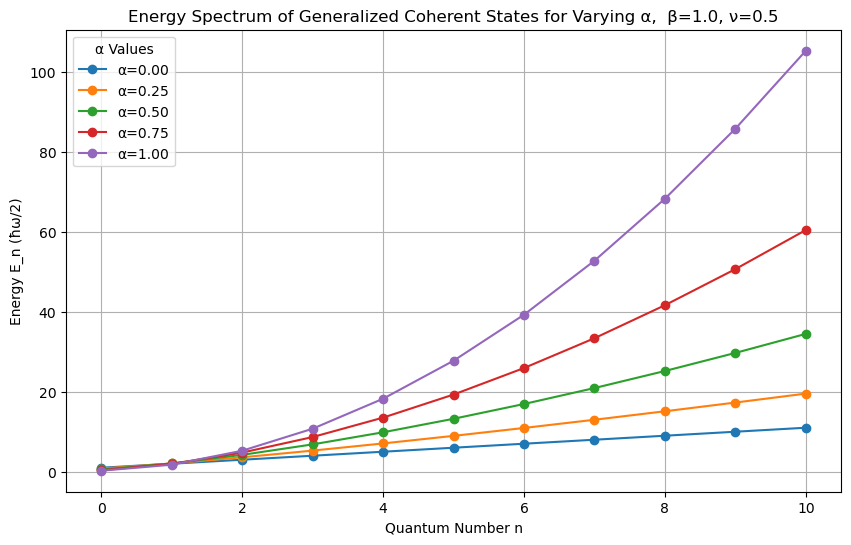}
        %\caption{}
        \label{fig:first}
    \end{subfigure}
    \hfill % Adjust spacing
    % Second image
    \begin{subfigure}[b]{0.45\linewidth}
        \centering
        \includegraphics[width=\linewidth]{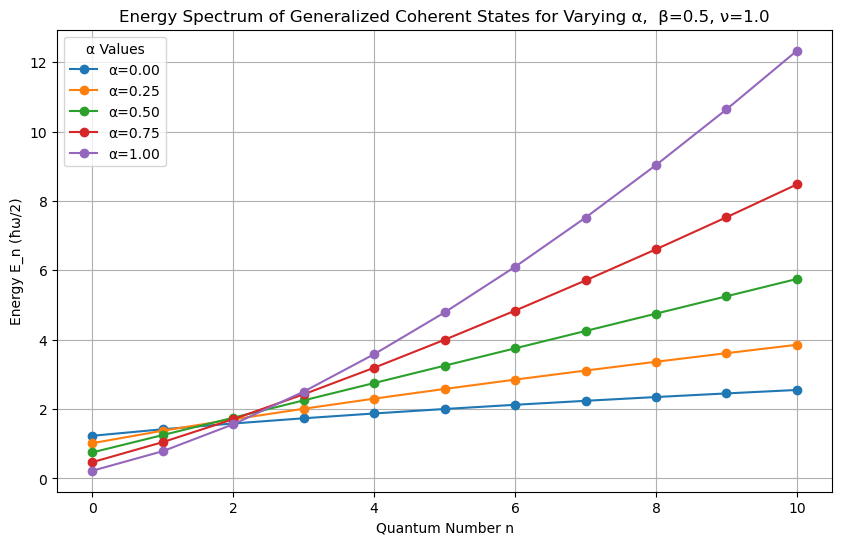}
       % \caption{}
        \label{fig:second}
    \end{subfigure}
    \caption{Energy spectra of generalised CS (\ref{energyev}) as a function of quantum numbers $n$ for varying deformation parameter $\alpha$. The plots illustrate the influence of the parameter $\alpha$ on the energy eigenvalues $E_n$, with fixed values of $\beta$ and $\nu$. Each curve corresponds to a specific value of $\alpha$, demonstrating the deformation effects in the spectral structure}
    \label{fig:energy spectra}
\end{figure}

$\{1,\hat A_{\alpha,\beta,\nu},\hat A^\dag_{\alpha,\beta,\nu},[\hat N]_{\alpha,\beta,\nu}\}$, where the deformed boson operators of annihilation $\hat A_{\alpha,\beta,\nu}$ and creation $\hat A^\dag_{\alpha,\beta,\nu}$ satisfy the deformed commutation rules:

\begin{equation}\label{commutation}
[\hat A_{\alpha,\beta,\nu}\,,\,\hat A^\dag_{\alpha,\beta,\nu}]=[\hat N+1]_{\alpha,\beta,\nu}-[\hat N]_{\alpha,\beta,\nu},
\end{equation}

%\begin{equation}
%[N_{\alpha,\beta,\nu},A_{\alpha,\beta,\nu}]=-A_{\alpha,\beta,\nu},\,\,\,\,\,\,\,\,\,\,[N_{\alpha,\beta,\nu},A^\dag_{\alpha,\beta,\nu}]=A^\dag_{\alpha,\beta,\nu},
%\end{equation}

where the { deformed} number operator $[\hat N]_{\alpha,\beta,\nu}$ is given by
 
\begin{equation}
[\hat N]_{\alpha,\beta,\nu}:=\hat A^\dag_{\alpha,\beta,\nu}\hat A_{\alpha,\beta,\nu}.
\end{equation}
{We note that this generalised perspective has also been explored in a different context by Solomon in \cite{Solomon94} and by Daskaloyannis and colleagues in \cite{Dask91,Dask92,DY92}}.
{We assumed that the number states $|n>$, elements of the Fock space, form an orthonormal basis of the deformed number operator $[\hat N]_{\alpha,\beta,\nu}$. (In the simplest case $\alpha=\nu=0,\,\beta=1$, the number operator $[\hat N]_{0,1,0}$ reduces to the standard operators of the Weyl-Heisenberg algebra.)} 
\begin{equation}
[\hat N]_{\alpha,\beta,\nu}|n>=[n]_{\alpha,\beta,\nu}|n>;
\end{equation}
which gives the following representation of the Fock space:
\begin{eqnarray}\label{creation-annihilation}
\hat A_{\alpha,\beta,\nu}|n>=\sqrt{[n]_{\alpha,\beta,\nu}}|n-1>;\\
\hat A^\dag_{\alpha,\beta,\nu}|n>=\sqrt{[n+1]_{\alpha,\beta,\nu}}|n+1>.
\end{eqnarray}
%The formula generates the eigenvectors of the number operator $N$
{However, we can use the deformed creation operator to construct all the other solutions for different $n$. Successive application of $\hat A^\dag_{\alpha,\beta,\nu}$ to the lowest state, corresponding to $n=0$, gives the normalized eigenstates:}
$$|n>=\frac{1}{\sqrt{[n]_{\alpha,\beta,\nu}!}}\left(\hat A^\dag_{\alpha,\beta,\nu}\right)^n|0>.$$
The number operator commutes with the following Hamiltonian
\begin{equation}\label{hamiltonian}
\hat H_{\alpha,\beta,\nu}=\frac{\hbar\omega}{2}\left(\hat A^\dag_{\alpha,\beta,\nu}\hat A_{\alpha,\beta,\nu}+\hat A_{\alpha,\beta,\nu}\hat A^\dag_{\alpha,\beta,\nu}\right);
\end{equation}
corresponding to the energy eigenvalues (see Figure \ref{fig:energy spectra})
\begin{equation}\label{energyev}
E_n^{\alpha,\beta,\nu}=\frac{\hbar\omega}{2}\left([n+1]_{\alpha,\beta,\nu}+[n]_{\alpha,\beta,\nu}\right);
\end{equation}
and the commutation rules between the Hamiltonian, the creation and the annihilation operators:

\begin{eqnarray}
\left[\hat H_{\alpha,\beta,\nu},\hat A^\dag_{\alpha,\beta,\nu}\right]&=&\frac{\hbar\omega}{2}\left([\hat N+1]_{\alpha,\beta,\nu}-[\hat N]_{\alpha,\beta,\nu}+1\right)\hat A^\dag_{\alpha,\beta,\nu};\\
\left[\hat H_{\alpha,\beta,\nu},\hat A_{\alpha,\beta,\nu}\right]&=&-\frac{\hbar\omega}{2}\left([\hat N+1]_{\alpha,\beta,\nu}-[\hat N]_{\alpha,\beta,\nu}+1\right)\hat A_{\alpha,\beta,\nu}.
\end{eqnarray}

%\begin{rem}	
%In case $\alpha=\nu=0$ and $\beta=1$ the Hamiltonian (\ref{hamiltonian}) becomes the Hamiltonian of the quantum harmonic %oscillator with discrete energy levels equally spaced $E_n=\hbar\omega\left(n+\frac{1}{2}\right)$.
%\end{rem}

\begin{defn}
Using $z^n\,\,\,(n=0,1,2...)$ as a basis function for this $\mathcal{W}_{\alpha,\beta,\nu}$-deformed boson algebra, we define the $\mathcal{W}_{\alpha,\beta,\nu}$-coherent states $|z;\alpha,\beta,\nu\rangle$ {(with $(\alpha\in [0,1],\beta\in (0,1])$ and $\nu\ge0$)} by the expression
\begin{equation}\label{CS}
|z;\alpha,\beta,\nu\rangle=[\mathcal{N}_{\alpha,\beta,\nu}(|z|^2)]^{-1/2}\sum_{n=0}^{\infty}\frac{z^n}{\sqrt{[n]_{\alpha,\beta,\nu}!}}|n\rangle,\,\,\,\,\,z\in\mathbb{C};
\end{equation}
{with the function $\mathcal{N}_{\alpha,\beta,\nu}(x)=\Gamma(1-\alpha+\nu)\mathcal{W}_{\alpha,\beta,\nu}(x)=\sum_{n=0}^{\infty}\frac{x^n}{[n]_{\alpha,\beta,\nu}!}$}, (see
 \cite{Droghei21} for more details about the $\mathcal{W}_{\alpha,\beta,\nu}$-function)
$$\mathcal{W}_{\alpha,\beta,\nu}(x)=\sum_{n=0}^\infty\prod_{i=1}^n\frac{\Gamma(\beta i+1-\alpha)}{\Gamma(\beta i+1)}\frac{x^{ n}}{\Gamma(\beta n+1-\alpha+\nu)}.$$

\begin{rem}
    {We define a fractional differential operator $_x\hat{D}_{\alpha,\beta,\nu}$ ( $\mathcal{W}_{\alpha,\beta,\nu}$-derivative) for real numbers:
\begin{equation}
_x\hat{D}_{\alpha,\beta,\nu}=x^{\alpha-\nu}\frac{d^\beta}{dx^\beta}x^\nu\frac{d^\alpha}{dx^\alpha};
\end{equation}
where $\frac{d^\gamma}{dx^\gamma}$ denotes the Caputo fractional derivative of order $\gamma$ . 
{ It acts as a deformed differential operator}
\begin{equation}
_x\hat{D}_{\alpha,\beta,\nu} \,\,x^{\beta n} = [n]_{\alpha,\beta,\nu}\,\,x^{\beta (n-1)}.
\end{equation}
Moreover, the $ \mathcal{W}_{\alpha,\beta,\nu}$-function is an eigenfunction of this operator. 
\begin{equation}
_x\hat{D}_{\alpha,\beta,\nu}\,\, \mathcal{N}_{\alpha,\beta,\nu}(\lambda x^\beta)=\lambda\,\, \mathcal{N}_{\alpha,\beta,\nu}(\lambda x^\beta).
\end{equation}
It plays the role of a deformed exponential function.
}
\end{rem}

{Using Eq.(\ref{creation-annihilation}), it is simple to demonstrate that the deformed $\mathcal{W}_{\alpha,\beta,\nu}$-coherent states are eigenstates of the deformed annihilation operator, $\hat A_{\alpha,\beta,\nu}|z;\alpha,\beta,\nu\rangle=z|z;\alpha,\beta,\nu\rangle$, which is a modification of the }{standard annihilation operator in the undeformed quantum harmonic oscillator. Furthermore, since $\hat A_{\alpha,\beta,\nu}\mathcal{N}_{\alpha,\beta,\nu}(z  \hat A^\dag_{\alpha,\beta,\nu})|0\rangle=z \mathcal{N}_{\alpha,\beta,\nu}(z  \hat A^\dag_{\alpha,\beta,\nu})|0\rangle$, we can use $\mathcal{N}_{\alpha,\beta,\nu}(x)$ to define analogues of coherent states as normalized eigenstates of the deformed annihilation operator
\begin{equation}
   |z;\alpha,\beta,\nu\rangle= \mathcal{N}_{\alpha,\beta,\nu}(|z|^2)^{-1/2}\mathcal{N}_{\alpha,\beta,\nu}(z   \hat A^\dag_{\alpha,\beta,\nu})|0\rangle.
\end{equation}
This observation allows us to place our work within the framework of deformed bosons.} 
\end{defn}

From the eq. (\ref{CS}), the probability of finding the state $|n>$ in the state ket $|z;\alpha,\beta,\nu>$  is equal to 

\begin{equation}
p_{\alpha,\beta,\nu}(n,z)=\frac{\Gamma(1-\alpha+\nu)\left(|z|^2\right)^n}{\left(\prod_{i=1}^n\frac{\Gamma(\beta i+1)}{\Gamma(\beta i+1-\alpha)}\right)\Gamma(\beta n+1-\alpha+\nu)\mathcal{W}_{\alpha,\beta,\nu}(|z|^2)}.
\end{equation}
It coincides with the Poisson distribution characterising the conventional CS for $\alpha=0,\, \beta=1, \,\nu=0$.

For two different complex numbers $z$ and $z'$, the states $|z\rangle$ and $|z'\rangle$ are, in general, not orthogonal and their overlap is given by $$\langle z|z'\rangle=\frac{\mathcal{N}_{\alpha,\beta,\nu}(z^*z')}{\sqrt{\mathcal{N}_{\alpha,\beta,\nu}(|z|^2)\cdot\mathcal{N}_{\alpha,\beta,\nu}(|z'|^2)}}$$ where $\mathcal{N}_{\alpha,\beta,\nu}(z^*z')=\sum_{n=0}^\infty\frac{\left(z^*z'\right)^n}{[n]_{\alpha,\beta,\nu}!}$.

\subsection{Deformed harmonic oscillator}
{In this section, we aim to investigate, in analogy to the {undeformed case}, the deformed boson algebra previously introduced about a deformed version of the quantum harmonic oscillator. To do this, we can define the position and momentum operators, $\hat q_{\alpha,\beta,\nu}$ and  $\hat p_{\alpha,\beta,\nu}$ respectively, for the generalised oscillator, concerning the creation and annihilation operators by:
\begin{equation}\label{q&p}
    \hat q_{\alpha,\beta,\nu}=\sqrt{\frac{\hbar}{2m\omega}}\left(\hat A_{\alpha,\beta,\nu}+\hat A^\dag_{\alpha,\beta,\nu}\right);\,\,\,\,\,\hat p_{\alpha,\beta,\nu}=-i\sqrt{\frac{\hbar m\omega}{2}}\left( \hat A_{\alpha,\beta,\nu}-\hat A^\dag_{\alpha,\beta,\nu}\right).
\end{equation}
Using the commutation rule (\ref{commutation}), we obtain
\begin{equation}
\left[\hat q_{\alpha,\beta,\nu}, \hat p_{\alpha,\beta,\nu}\right]= i\hbar[\hat A_{\alpha,\beta,\nu}\,,\,\hat A^\dag_{\alpha,\beta,\nu}]=i\hbar\left([\hat N+1]_{\alpha,\beta,\nu}-[\hat N]_{\alpha,\beta,\nu}\right);
\end{equation}
and the deformed Heisenberg's equations of motion:
\begin{eqnarray}
\dot{\hat q}_{\alpha,\beta,\nu} &=& \left([n+1]_{\alpha,\beta,\nu}-[n]_{\alpha,\beta,\nu}+1\right)\frac{\hat p_{\alpha,\beta,\nu}}{2m};\nonumber\\
\dot{\hat p}_{\alpha,\beta,\nu} &=& -\left([n+1]_{\alpha,\beta,\nu}-[n]_{\alpha,\beta,\nu}+1\right)\frac{m\omega^2 \hat q_{\alpha,\beta,\nu}}{2}.
\end{eqnarray}
}
{ In analogy with the classical undeformed case, we define the deformed creation and annihilation operators as follows:
\begin{equation}\label{creation-annihilation-differential}
\hat A^\dag_{\alpha,\beta,\nu}:=\left(\sqrt{\frac{m\omega}{2\hbar}}x^\beta-\sqrt{\frac{\hbar}{2 m\omega}}\,_x\hat{D}_{\alpha,\beta,\nu}\right),\,\,\,\,\,\hat A_{\alpha,\beta,\nu}:=\left(\sqrt{\frac{m\omega}{2\hbar}}x^\beta+\sqrt{\frac{\hbar}{2 m\omega}}\,_x\hat{D}_{\alpha,\beta,\nu}\right);
\end{equation}
where, in the undeformed case, the momentum operator is given by  $\hat p= -i\hbar\frac{d}{dx}$.
By imposing the ground state condition $\hat A_{\alpha,\beta,\nu}|0\rangle=0$, we obtain the analytical expression of the eigenfunction by solving the fractional differential equation: $$\left(\sqrt{\frac{m\omega}{2\hbar}}x^\beta+\sqrt{\frac{\hbar}{2 m\omega}}\,_x\hat{D}_{\alpha,\beta,\nu}\right)\langle x|0\rangle=0.$$
It is straightforward to show that the normalized ground state wave function takes the form:
\begin{equation}\label{groundstate}
    \langle x|0\rangle=\left(\frac{m\omega}{\pi\hbar}\right)^{\frac{1}{4}}\sum_{n=0}^\infty\frac{\left(-\frac{m\omega}{\hbar}\right)^n}{[2n]_{\alpha,\beta,\nu}!!}x^{2\beta n},
\end{equation}

where $[n]_{\alpha,\beta,\nu}!!$ denotes the generalised double factorial associated with (\ref{genfact}).
To determine the eigenfunctions of the excited states, we recursively apply the creation operator to the ground state $ \langle x|n\rangle=\frac{1}{\sqrt{[n]_{\alpha,\beta,\nu}!}}\left(\hat A^\dag_{\alpha,\beta,\nu}\right)^n|0\rangle$.  For instance, the first excited state is given by
\begin{equation}\label{firststate}
    \langle x|1\rangle=\sqrt{\frac{m\omega}{2\hbar}}\frac{2\,x^\beta}{\sqrt{[1]_{\alpha,\beta,\nu}!}}\langle x|0\rangle.
\end{equation}
}

\begin{figure}[ht]
    \centering
    \includegraphics[width=0.95\linewidth]{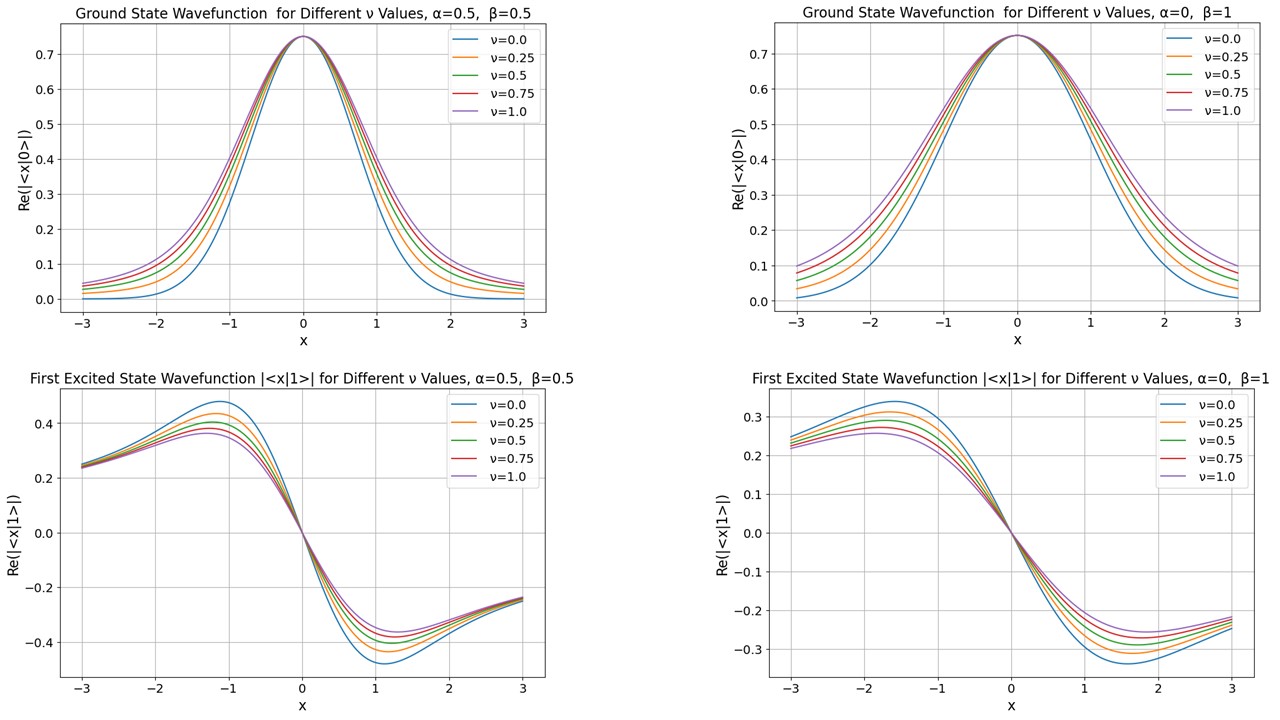}
    \caption{Ground states $ \langle z|0\rangle$  and First excited states  $\langle z|1\rangle $ for different values of the deformation parameters ($\alpha, \beta, \nu$).  We set: $m=\omega=\hbar=1$}
    \label{fig:eigenfunctions}
\end{figure}

\begin{rem}
    Obviously, for $\alpha=\nu=0$ and $\beta=1$, we obtain the eigenfunctions of the classical harmonic oscillator.
\end{rem}
{Figure~\ref{fig:eigenfunctions} illustrates the behavior of the wavefunctions of the ground state \( \langle x|0\rangle \) and the first excited state \( \langle x|1\rangle \) for different values of the deformation parameter \( \nu \), in two sets of distinct parameters (\(\alpha = 0, \beta = 1\) and \(\alpha = 0.5, \beta = 0.5\)). The plots reveal how variations in \( \nu \) influence the symmetry, amplitude, and spatial distribution of the wavefunctions, highlighting the effects of deformation on both ground and excited states.}

\subsection{Continuity and completeness of $\mathcal{W}_{\alpha,\beta,\nu}$-coherent states}
{ Although the term coherent state applies to a broad class of objects, all definitions share two fundamental properties, continuity and completeness (resolution of unity), which can be considered the minimal requirements for a set of states to be classified as coherent. The two postulates about CS were proposed in substantially this form nearly sixty years ago in \cite{Klauder63} }. FIRST CONDITION: the states (\ref{CS}) are coherent if they are continuous in the label $z$. {Continuity retains its standard meaning: for any convergent sequence of labels set $\mathcal{Z}$ such that $z'-z\rightarrow 0$, it follows that $\||z'\rangle-|z\rangle\|\rightarrow0$. As usual, the vector norm is defined as $ | |\psi\rangle | =\sqrt{\langle \psi | \psi \rangle}$, something that is always positive except in trivial cases $|\psi\rangle = 0 $. Furthermore, we assume that $|z\rangle \neq 0 \,\, \text{for all} \,\, z \in \mathcal{Z}$  }. This condition follows from the joint continuity of the reproducing kernel $\mathcal{K}(z;z')\equiv\langle z|z'\rangle=\langle z'|z\rangle^*$, in fact:
\begin{equation}\label{continuity}
\||z'\rangle-|z\rangle\|^2=\langle z'|z'\rangle+\langle z|z\rangle-2\Re{\langle z'|z\rangle}=2(1-\Re{\langle z'|z\rangle}), 
\end{equation}
 and is easily satisfied in practice, see \cite{KS85}.

SECOND CONDITION: the conditions of completeness, and hence the resolution of unity (\ref{completeness}), imposes for $x\equiv|z|^2$ the equation
\begin{equation}\label{Stiel}
\int_0^\infty x^n\left[\pi\frac{U(x)}{\mathcal{N}_{\alpha,\beta,\nu}(x)}\right]dx=[n]_{\alpha,\beta,\nu}!.
\end{equation}
(see \cite{KPS2001} for more details). {The quantities $[n]_{\alpha,\beta,\nu}!$ are then the power of the unknown positive function $\tilde{U}(x)=\pi\frac{U(x)}{\mathcal{N}_{\alpha,\beta,\nu}(x)}$ and the problem stated in equation (\ref{Stiel}) is the Stieltjes moment problems.} It is important to note that not all deformed algebras lead to CS within the construction framework, as the moment problem (\ref{Stiel}) does not always have a solution, see \cite{A1965},\cite{ST1943}. { As we explain in the following subparagraph, the deformed algebra we investigated leads to CS only if the deformation parameters satisfy Carleman's conditions, namely $\alpha+\beta\le 2$ and $\nu>\alpha-1$}

We know that strictly positive determinants of the Hankel-Hadamard matrices are necessary and sufficient for the weight function to exist \cite{SPS99},\cite{A1965}.

The not-trivial Stieltjes problem can be tackled using a Mellin and inverse Mellin transforms approach, extending the natural value $n$ to the complex values $s$ and rewriting (\ref{Stiel}) as
\begin{equation}
\int_0^\infty x^{s-1}\tilde{U}(x)dx=[s-1]_{\alpha,\beta,\nu}!.
\end{equation}

\begin{rem}
If $\alpha=0$, we obtain $[n]_{0,\beta,\nu}!=\frac{\Gamma(\beta n+1+\nu)}{\Gamma(1+\nu)}$ which is the generalised factorial performed by adopting the Mittag-Leffler function.
If $\alpha=1$, $[n]_{1,\beta,\nu}!=\beta^n n!\frac{\Gamma(\beta n+\nu)}{\Gamma(\nu)}$ we have the generalised factorial in the case of the classical Wright function.
\end{rem}

%---------------------------------------------------------------------------------------------------------------------------------------------------------------------------------------
\subsubsection{Resolution Of Unity}
From the definition of finite moments:

\[
[n]_{\alpha,\beta,\nu}! = \left(\prod_{i=1}^{n} \frac{\Gamma(\beta i + 1)}{\Gamma(\beta i + 1 - \alpha)}\right) \cdot \frac{\Gamma(\beta n + 1 - \alpha + \nu)}{\Gamma(1 - \alpha + \nu)}.
\]
We aim to identify the conditions necessary for the existence of a measure corresponding to these moments.
The moment sequence presented involves products of ratios of Gamma functions. To analyze the existence of a measure linked to these moments, it is essential to ensure that the moment sequence behaves appropriately, particularly in growth and positivity, while verifying the applicability of established conditions. Carleman's condition offers a sufficient criterion for the determinacy of the moment problem.
In the Appendix A, we provide a detailed analysis of Carleman's condition concerning our finite moments $[n]_{\alpha,\beta,\nu}!$ 
Specifically, we find that the condition $\frac{\alpha+\beta}{2}\le 1$ is sufficient to ensure that Carleman's condition holds, thus confirming the uniqueness of the associated measure.\cite{A1965}

%---------------------------------------------------------------------------------------------------------------------------------------------------------------------------------------

\subsection{Positive weight function in the Wright coherent states}%\label{sec:2}

The family of CS using the Mittag-Leffler functions has been deeply investigated by Sixdeniers and collaborators in \cite{SPS99}.
The Wright function case has been analysed in \cite{GGM19} and \cite{GM23}.
In the following, we solve the Stieltjes problem in the case of Wright generalised factorial ($\alpha=1$), considering the following equation as the Melling transform for complex variable $s$, obtaining the same result given by Giraldi and Mainardi in \cite{GM23}

\begin{equation}\label{Mell}
\int_0^\infty x^{s-1}\tilde{U}_{1,\beta,\nu}(x)dx=[s-1]_{1,\beta,\nu}!=\beta^{s-1}\cdot \frac{\Gamma(s)\Gamma(\beta(s-1)+\nu)}{\Gamma(\nu)}.
\end{equation}
To obtain $\tilde{U}_{1,\beta,\nu}(x)$, we have to carry out an inverse Mellin transform on (\ref{Mell}).

\begin{equation}\label{inMell}
\tilde{U}_{1,\beta,\nu}(x)=\frac{1}{\beta\Gamma(\nu)}\frac{1}{2\pi\imath}\int_{\mathcal{L}}\left(\frac{x}{\beta}\right)^{-s}\Gamma(s)\Gamma(\beta(s-1)+\nu)ds=\frac{1}{\beta\Gamma(\nu)}H_{0,2}^{2,0}\left[\frac{x}{\beta}\left| \begin{array}{ll}
         -\\
        (0,1)\,\,,\,\,(\nu-\beta,\beta).\end{array} \right.\right]
\end{equation}

that is the Fox-H function 
\begin{equation*}
    H_{p,q}^{m,n}\left[z\left| 
    \begin{array}{ll}
         (a_1,A_1),...,(a_p,A_p)\\
        (b_1,B_1),...,(b_q,B_q)\end{array} \right.\right]=\frac{1}{2\pi\imath}\int_{\mathcal{L}}\frac{\prod_{j=1}^m\Gamma(b_j+B_js)\prod_{j=1}^n\Gamma(1-a_j-A_js)}{\prod_{j=m+1}^q\Gamma(1-b_j-B_js)\prod_{j=n+1}^p\Gamma(a_j+A_js)}z^{-s}ds.
\end{equation*}

Refer to \cite{MSH09},\cite{K2004} for details of the general theory and application of the H-function.
Considering  the relations (2.9.31)and (2.9.32) in \cite{K2004}, we can rewrite the equation (\ref{inMell}) in the following integral form

\begin{equation}
\tilde{U}_{1,\beta,\nu}(x)=\frac{1}{\beta^2\Gamma(\nu)}\int_0^\infty
t^{\frac{\nu}{\beta}-2}exp\left(-t^{\frac{1}{\beta}}-\frac{x}{\beta t}\right)dt.
\end{equation}

\subsubsection{Case $\alpha=\beta$}
\begin{equation}\label{Mell2}
\int_0^\infty x^{s-1}\tilde{U}_{\beta,\beta,\nu}(x)dx=[s-1]_{\beta,\beta,\nu}!= \frac{\Gamma(\beta s+1-\beta)\Gamma(\beta(s-2)+1+\nu)}{\Gamma(1-\beta+\nu)}.
\end{equation}
To obtain $\tilde{U}_{\beta,\beta,\nu}(x)$, we have to perform an inverse Mellin transform on (\ref{Mell2}).

\begin{eqnarray}\label{inMell2}
\tilde{U}_{\beta,\beta,\nu}(x)=\frac{1}{\Gamma(1-\beta+\nu)}\frac{1}{2\pi\imath}\int_{\mathcal{L}}x^{-s}\Gamma(\beta s+1-\beta)\Gamma(\beta(s-2)+1+\nu)ds\\
=\frac{1}{\Gamma(1-\beta+\nu)}H_{0,2}^{2,0}\left[\frac{x}{\beta}\left| \begin{array}{ll}
         -\\
        (1-\beta,\beta)\,\,,\,\,(1-2\beta+\nu,\beta).\end{array} \right.\right]
\end{eqnarray}

\subsubsection{Case $\alpha=1-\beta$}
\begin{equation}\label{Mell3}
\int_0^\infty x^{s-1}\tilde{U}_{1-\beta,\beta,\nu}(x)dx=[s-1]_{1-\beta,\beta,\nu}!= \frac{\beta^{s-1}\Gamma(\beta)\Gamma(s)\Gamma(\beta s+\nu)}{\Gamma(\beta+\nu)\Gamma(\beta s)}.
\end{equation}
To obtain $\tilde{U}_{1-\beta,\beta,\nu}(x)$, we have to perform an inverse Mellin transform on (\ref{Mell3}).

\begin{eqnarray}\label{inMell3}
\tilde{U}_{1-\beta,\beta,\nu}(x)=\frac{\Gamma(\beta)}{\beta\Gamma(\beta+\nu)}\frac{1}{2\pi\imath}\int_{\mathcal{L}}\left(\frac{x}{\beta}\right)^{-s} \frac{\Gamma(s)\Gamma(\beta s+\nu)}{\Gamma(\beta s)}ds\\
=\frac{\Gamma(\beta)}{\beta\Gamma(\beta+\nu)}H_{1,2}^{2,0}\left[\frac{x}{\beta}\left| \begin{array}{ll}
         (0,\beta)\\
        (0,1)\,\,,\,\,(\nu,\beta).\end{array} \right.\right]\\
=\frac{\Gamma(\beta)}{\beta^2\Gamma(\beta+\nu)\Gamma(-\nu)}\int_1^\infty (t^{\frac{1}{\beta}}-1)^{-\nu-1}t^{\frac{1}{\beta}-1}e^{-\frac{x}{\beta}t}dt.
\end{eqnarray}
Other two subcases $[s]_{0,2r,0}!=\Gamma(2rs+1)$ and $[s]_{\beta,r,\beta}!=(\Gamma(rs+1))^2$ have been investigated  by Penson et al. in \cite{PBDHS2010}

\section{Quantum Fluctuations of Quadrature}\label{sec:3}
\setcounter{section}{3} \setcounter{equation}{0}
Characterising states through their quadrature operators in quantum mechanics provides crucial insights into the inherent uncertainties that govern their behaviour. The quadrature operators denoted here as \(q_{\alpha,\beta,\nu}\) and \(p_{\alpha,\beta,\nu}\) defined in (\ref{q&p}), are fundamental in analysing the properties of generalised CS, particularly with respect to their fluctuations. These operators represent position- and momentum-like variables and are defined as the corresponding annihilation and creation operators.

We rewrite the explicit forms of the quadrature operators:

\begin{equation}
\hat q_{\alpha,\beta,\nu}=\sqrt{\frac{\hbar}{2m\omega}}\left(\hat A_{\alpha,\beta,\nu}+ \hat A^\dag_{\alpha,\beta,\nu}\right), \quad \hat p_{\alpha,\beta,\nu}=i\sqrt{\frac{\hbar m\omega}{2}}\left(\hat A^\dag_{\alpha,\beta,\nu}-\hat A_{\alpha,\beta,\nu}\right).
\end{equation}

The expectation values of these operators are derived in terms of the state parameters, providing a foundation for calculating their variances which quantify the quantum fluctuations. Specifically, we find:

\begin{align}
\langle \hat q_{\alpha,\beta,\nu} \rangle &= \sqrt{\frac{\hbar}{2m\omega}}\left(\langle \hat A_{\alpha,\beta,\nu} \rangle + \langle \hat A^\dag_{\alpha,\beta,\nu} \rangle\right), \\
\langle \hat p_{\alpha,\beta,\nu} \rangle &= i\sqrt{\frac{\hbar m\omega}{2}}\left(\langle \hat A^\dag_{\alpha,\beta,\nu} \rangle - \langle \hat A_{\alpha,\beta,\nu} \rangle\right).
\end{align}

The fluctuations in the quadratures are quantified by the variances \(\Delta \hat q_{\alpha,\beta,\nu}\) and \(\Delta \hat p_{\alpha,\beta,\nu}\), which encapsulate the spread of the measurements around the mean values. In particular, we derive the relations for the variances as follows:

\begin{align}
\Delta \hat q^2_{\alpha,\beta,\nu} &= \frac{\hbar}{2m\omega}\{[n+1]_{\alpha,\beta,\nu}-[n]_{\alpha,\beta,\nu}\}, \\
\Delta \hat p^2_{\alpha,\beta,\nu} &= \frac{\hbar m\omega}{2}\{[n+1]_{\alpha,\beta,\nu}-[n]_{\alpha,\beta,\nu}\}.
\end{align}
In the vacuum state
\begin{align}
(\Delta \hat q_{\alpha,\beta,\nu})_0 &= \sqrt{\frac{\hbar}{2m\omega}[1]_{\alpha,\beta,\nu}}, \\
(\Delta \hat p_{\alpha,\beta,\nu})_0 &=\sqrt{ \frac{\hbar m\omega}{2}[1]_{\alpha,\beta,\nu}}.
\end{align}
These results culminate in the expression of the product of the uncertainties, which adheres to the Heisenberg uncertainty principle:

\begin{equation}
\Delta q_{\alpha,\beta,\nu} \Delta p_{\alpha,\beta,\nu} = \frac{\hbar}{2}\{[n+1]_{\alpha,\beta,\nu}-[n]_{\alpha,\beta,\nu}\}.
\end{equation}
Finally, we express the uncertainty relation for our CS in terms of the generalised gamma functions, emphasising the dependence on the state parameters:\\

\scriptsize
\begin{equation}
\Delta \hat q_{\alpha,\beta,\nu}\Delta \hat p_{\alpha,\beta,\nu}=\frac{\hbar}{2}\left[\frac{\Gamma(\beta(n+1)+1)}{\Gamma(\beta(n+1)+1-\alpha)}\frac{\Gamma(\beta(n+1)+1-\alpha+\nu)}{\Gamma(\beta n+1-\alpha+\nu)} - \frac{\Gamma(\beta n+1)}{\Gamma(\beta n+1-\alpha)}\frac{\Gamma(\beta n+1-\alpha+\nu)}{\Gamma(\beta (n-1)+1-\alpha+\nu)}\right]\\
\quad n\ge 1.
\end{equation}
\normalsize

\begin{figure}[ht]
    \centering
    % First image
    \begin{subfigure}[b]{0.45\linewidth}
        \centering
        \includegraphics[width=\linewidth]{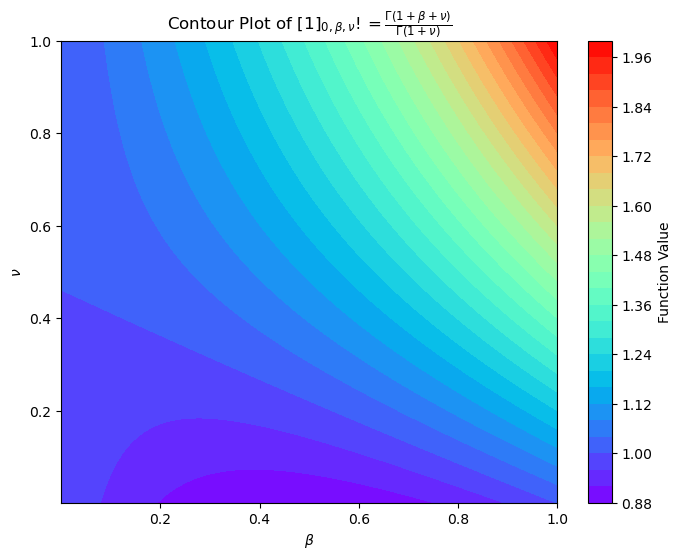}
        \caption{Mittag-Leffler case: $\alpha = 0$}
        \label{fig:first3}
    \end{subfigure}
    \hfill % Adjust spacing
    % Second image
    \begin{subfigure}[b]{0.45\linewidth}
        \centering
        \includegraphics[width=\linewidth]{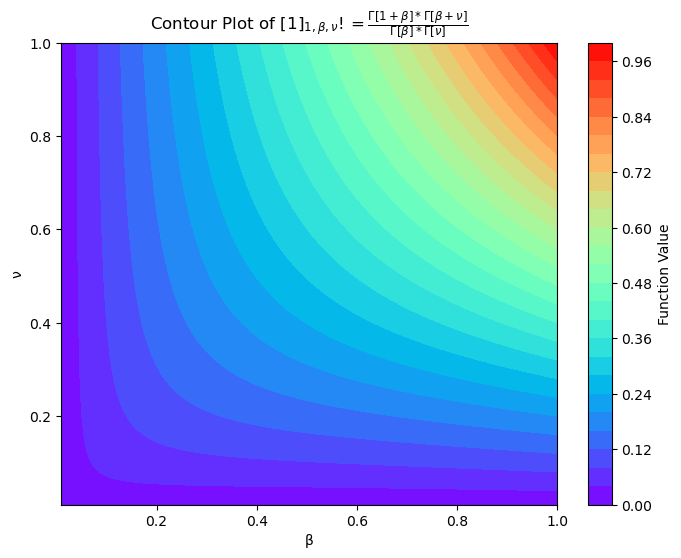}
        \caption{Wright function case: $\alpha = 1$}
        \label{fig:second3}
    \end{subfigure}
    \caption{Visualisation of quantum fluctuations in quadrature uncertainties of the vacuum state for the generalised deformed states, expressed in units of $\hbar/2$, illustrating the dependence of the uncertainty relation (Equation 3.10) on the parameters $\alpha$, $\beta$, and $\nu$.}
    \label{fig:uncertainty_relation}
\end{figure}

\begin{figure}[!]
    \centering
    \includegraphics[width=0.5\linewidth]{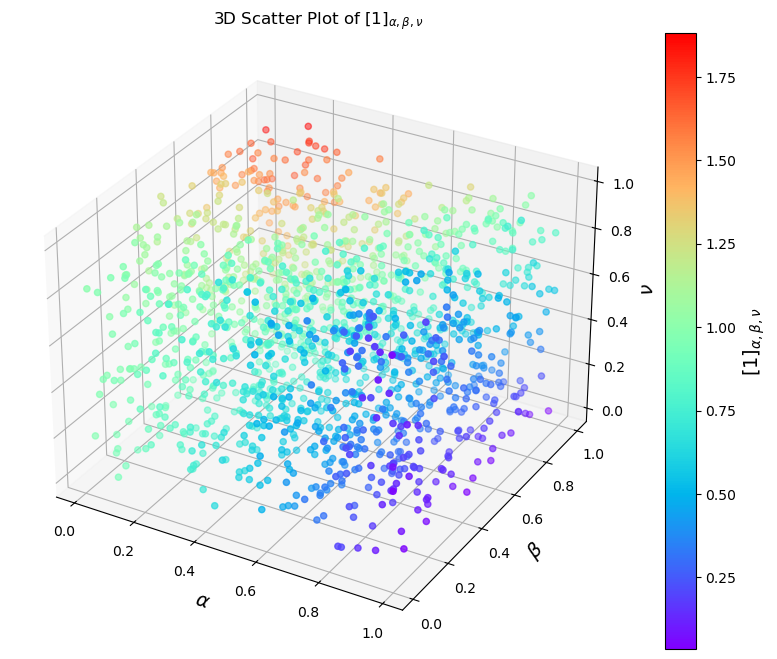}
    \caption{3D scatter plot illustrating the dependency of the gamma ratio 
$[1]_{\alpha,\beta,\nu}=\frac{\Gamma(\beta +1)}{\Gamma(\beta +1-\alpha)} \frac{\Gamma(\beta +1-\alpha+\nu)}{\Gamma(1-\alpha+\nu)}$ 
on the parameters $\alpha$, $\beta$, and $\nu$ within the range $(0, 1)$. The colour scale represents the computed gamma ratio values, highlighting variations across the parameter space.}
    \label{fig:3dscatter}
\end{figure}

\vspace{10pt}
This analysis not only elucidates the quantum fluctuations present in generalised CS but also underscores the interplay between the state parameters and the fundamental limits imposed by quantum mechanics. Through this framework, we gain valuable insights into the nature of quantum noise and its implications for quantum technologies. As shown in Figure~\ref{fig:uncertainty_relation},~\ref{fig:3dscatter}, the product of the uncertainties of the vacuum state 

\begin{equation}
\left(\Delta \hat q_{\alpha,\beta,\nu}\right)_0\left( \Delta \hat p_{\alpha,\beta,\nu}\right)_0 = \frac{\hbar}{2}[1]_{\alpha,\beta,\nu}= \frac{\hbar}{2}\frac{\Gamma(\beta +1)}{\Gamma(\beta +1-\alpha)}\frac{\Gamma(\beta +1-\alpha+\nu)}{\Gamma(1-\alpha+\nu)};
\end{equation}
depends on the parameters $\alpha$, $\beta$, and $\nu$. It shows that in the case $\alpha=0$ (Figure~\ref{fig:first3}), the uncertainty relation of the vacuum state goes from a value of $\frac{\hbar}{2}$ to $\hbar$ representing a higher level of quantum noise or fluctuations than the standard coherent state.

\section{Mandel parameter}\label{sec:4}

\setcounter{section}{4} \setcounter{equation}{0}
\begin{figure}[!b]
    \centering
    \begin{minipage}[b]{0.45\linewidth}
        \centering
        \includegraphics[width=\linewidth]{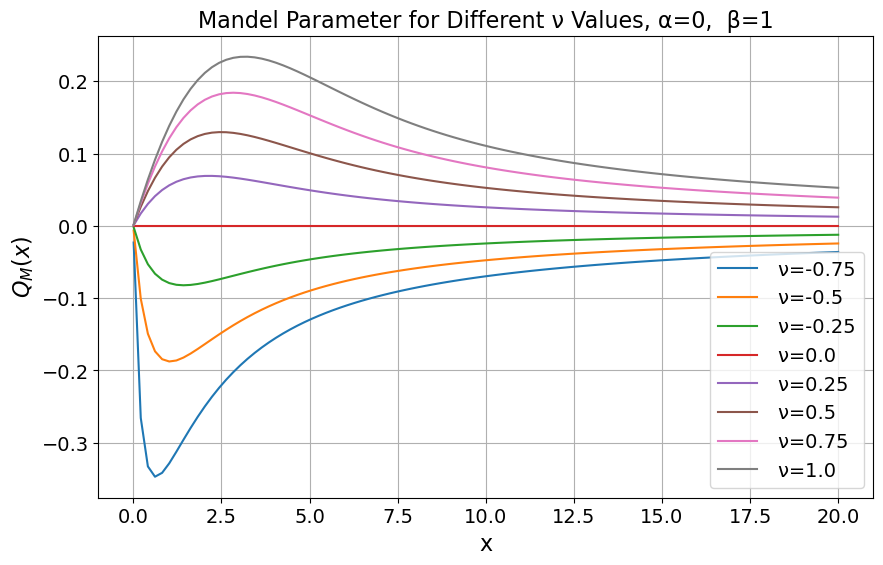}
        \caption{Plot of the Mandel parameter $Q_M(x)$ as a function of $x=|z|^2$, in the case of  $\alpha =0, \beta=1$ for different values of $\nu$. Mittag-Leffler function $E_{1,\nu+1}(x)$.}
    \label{fig:MandelML}
    \end{minipage}
    \hspace{0.05\linewidth}
    \begin{minipage}[b]{0.45\linewidth}
        \centering
        \includegraphics[width=\linewidth]{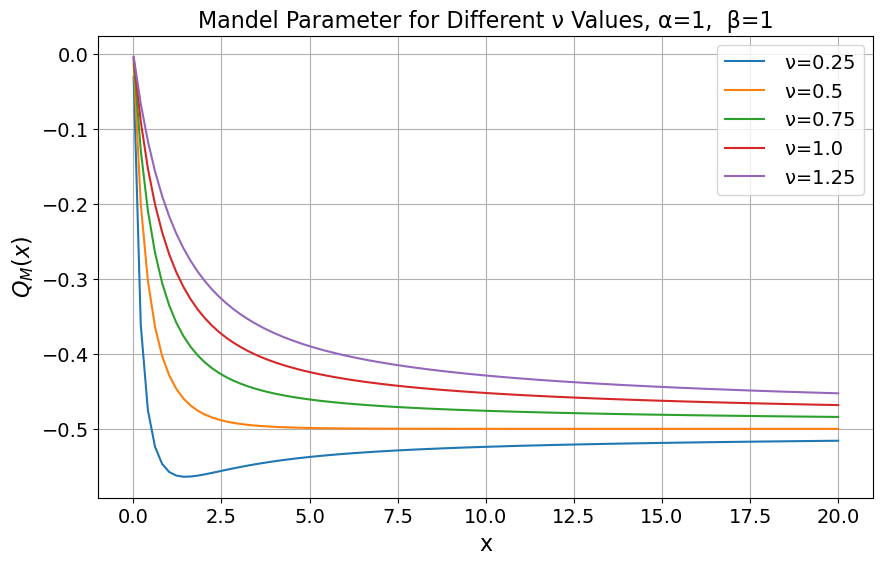}
        \caption{Plot of the Mandel parameter $Q_M(x)$ as a function of $x=|z|^2$, in the case of  $\alpha =1, \beta=1,$ for different values of $\nu$. Wright function $W_{1,\nu}(x)$}
     \label{fig:MandelWright}
    \end{minipage}
\end{figure}
While Glauber CS conventionally describes an ideal laser's states, real lasers do not strictly conform to this model. In particular, the photon number statistics of real lasers deviate from a Poissonian distribution, often due to various nonlinear interactions that lead to distinct departures from the ideal case. Recently, deformations of the commutation rules of boson operators have been proposed to model physical systems that derive from these idealised behaviours \cite{KS1994}, \cite{Solomon94}. The real laser problem was addressed in this phenomenological context, demonstrating that CS of deformed boson operators provides a more accurate model for non-ideal lasers, particularly concerning photon number statistics. A Poisson distribution is characterised by the variance of the deformed number operator $[\hat N]_{\alpha,\beta,\nu}$ being equal to its average. The deviation from the Poisson distribution can be measured with the Mandel parameter $Q_M(x)$, $(x=|z|^2)$ \cite{MW95}

\begin{equation}\label{mandel1}
Q_M=\frac{<[\hat N]_{\alpha,\beta,\nu}^2>-<[\hat N]_{\alpha,\beta,\nu}>^2}{<[\hat N]_{\alpha,\beta,\nu}>}-1.
\end{equation}

Using the relation between various expectation values of polynomial Hermitian operators and the derivative of $\mathcal{N}_{\alpha,\beta,\nu}(|z|^2)$, see also ref. \cite{KPS2001}

\begin{equation}\label{power}
<(\hat A^\dag_{\alpha,\beta,\nu})^r\hat A^r_{\alpha,\beta,\nu}>=\frac{|z|^{2r}}{\mathcal{N}_{\alpha,\beta,\nu}(|z|^2)}\left(\frac{d}{d|z|^2}\right)^r\mathcal{N}_{\alpha,\beta,\nu}(|z|^2),\,\,\,\,\,r=0,1,2,...\,.
\end{equation}

We have the following expectation values computed over the generalised CS $|z;\alpha,\beta,\nu>$:
\begin{equation}
<\hat A^\dag_{\alpha,\beta,\nu}\hat A_{\alpha,\beta,\nu}>=\frac{|z|^2}{\mathcal{W}_{\alpha,\beta,\nu}(|z|^2)}\sum_{n=0}^\infty\prod_{i=1}^{n+1}\frac{\Gamma(\beta i+1-\alpha)}{\Gamma(\beta i+1)}\frac{(n+1)|z|^{2 n}}{\Gamma(\beta (n+1)+1-\alpha+\nu)};
\end{equation}

\begin{equation}
<(\hat A^\dag_{\alpha,\beta,\nu})^2\hat A_{\alpha,\beta,\nu}^2>=\frac{|z|^4}{\mathcal{W}_{\alpha,\beta,\nu}(|z|^2)}\sum_{n=0}^\infty\prod_{i=1}^{n+2}\frac{\Gamma(\beta i+1-\alpha)}{\Gamma(\beta i+1)}\frac{(n+1)(n+2)|z|^{2 n}}{\Gamma(\beta (n+2)+1-\alpha+\nu)}.
\end{equation}
For the $\mathcal{W}_{\alpha,\beta,\nu}$-coherent states (\ref{CS}), the corresponding Mandel $Q_z$ parameter can be evaluated via the expectation values which are listed above and via the expression below
\begin{equation}
    Q_z=\frac{<(\hat A^\dag_{\alpha,\beta,\nu})^2\hat A_{\alpha,\beta,\nu}^2>-\left(<\hat A^\dag_{\alpha,\beta,\nu}\hat A_{\alpha,\beta,\nu}>\right)^2}{<\hat A^\dag_{\alpha,\beta,\nu}\hat A_{\alpha,\beta,\nu}>}.
\end{equation}
Negative values of the Mandel factor ($-1<Q_z<0$) (see Figures \ref{fig:MandelML} and \ref{fig:MandelWright}) indicate the nonclassical nature of the states, revealing the sub-Poissonian photon number statistic, a phenomenon with a nonclassical analogue. { Additionally, Figure \ref{fig:MandelML} shows, for positive values of  $0<\nu<1$, a Mandel's parameter $Q_z>0$, indicating a super-Poissonian statistic, which is typical of thermal states. For $\nu=0$,  $Q_z=0$, the field exhibits a Poissonian distribution (typical of the coherent light).} See \cite{MW95} for further details.

%\begin{figure}[th]
%    \centering
 %   \includegraphics[width=0.5\linewidth]{MandelML.pdf}
  %  \caption{Plot of the Mandel parameter $Q_M(x)$ as a function of $x=|z|^2$, in the case of  $\alpha =0, \beta=1,\,\nu=-1/2$. Mittag-Leffler %function $E_{1,\frac{1}{2}}(x)$}.
   % \label{fig:enter-label}
%\end{figure}
%\begin{figure}
 %   \centering
  %  \includegraphics[width=0.5\linewidth]{MandelWright.pdf}
   % \caption{Plot of the Mandel parameter $Q_M(x)$ as a function of $x=|z|^2$, in the case of  $\alpha =1, \beta=1,\,\nu=1/2$. Wright function $W_{1,\frac{1}{2}}(x)$.}
   % \label{fig:enter-label}
%\end{figure}

%\begin{figure}[h]
%%\centering
%\includegraphics[width=0.3\linewidth]{/Users/riccardodroghei/Documents/Droghei-Garra/Coherent-states/MandelWright.pdf}
%\caption{Plot of the Mandel parameter $Q_M(x)$ as a function of $x=|z|^2$, in the case of  $\alpha =1, \beta=1,\,\nu=1/2$.}
%\label{fig:bytepost}
%\end{figure}

%---------------------------------------------------------------------------------------------------------------------------------------------------------------------------------------
\section{Conclusion}
{
In this study, we introduced a novel class of coherent states, the $\mathcal{W}_{\alpha, \beta, \nu}$-coherent states, which extend the conventional framework of quantum optics and quantum mechanics by incorporating advanced mathematical structures from fractional calculus. The deformed boson algebra underlying these states is based on a generalized factorial function, allowing for a more flexible treatment of quantum systems where traditional commutation relations are insufficient.

Through this framework, we demonstrated that these deformed coherent states satisfy fundamental quantum properties such as continuity and completeness. We resolved the unity condition via a Stieltjes moment problem, which we addressed using Mellin and inverse Mellin transforms. Our analysis of quantum fluctuations in quadrature operators revealed a direct dependence on the deformation parameters $(\alpha, \beta, \nu)$, affecting position and momentum uncertainties. The evaluation of the Mandel parameter further underscored the non-classical properties of these states, highlighting sub-Poissonian photon statistics in specific parameter regimes.
The asymptotic analysis of the generalised factorials, as shown in Appendix A, further validated the consistency of our results, satisfying Carleman’s condition under specific parameter constraints.

By integrating fractional differential operators into the quantum framework, we have opened new possibilities for modelling quantum systems that exhibit memory effects, non-locality, and anomalous diffusion. These generalizations are not merely mathematical curiosities but hold potential experimental relevance, particularly in the study of fractional quantum mechanics, generalized uncertainty principles, and quantum optics. Future research may focus on practical implementations of these states in quantum information processing, non-classical light generation, and fractional field theories, as well as exploring other algebraic structures that emerge naturally in this extended quantum framework.

Despite the progress made, two key points remain open for future investigation. First, the physical nature and interpretation of the deformation parameters $(\alpha, \beta, \nu)$ require further clarification to fully understand their role in quantum mechanical formulations. Second, while the measure ensuring the completeness relation is unique up to a unitary transformation, the implications of non-unique measures and their relation to the choice of observables for different observers warrant deeper exploration. Future investigations will be crucial in addressing these open questions, leading to a more comprehensive understanding of the quantum mechanical framework developed in this study.

}

\section{Appendix A}

Let's analyse the asymptotic behaviour of the generalised factorial
$$
[n]_{\alpha,\beta,\nu} ! = \left( \prod_{i=1}^{n} \frac{\Gamma(\beta i + 1)}{\Gamma(\beta i + 1 - \alpha)} \right) \frac{\Gamma(\beta n + 1 - \alpha + \nu)}{\Gamma(1 - \alpha + \nu)}.
$$
We can break this generalised factorial into two parts:
The product term:
$$
P(n) = \prod_{i=1}^{n} \frac{\Gamma(\beta i + 1 )}{\Gamma(\beta i + 1- \alpha)}.
$$
The final term:
$$
F(n) = \frac{\Gamma(\beta n + 1 - \alpha + \nu)}{\Gamma(1 - \alpha + \nu)}.
$$
We'll analyse the asymptotic behaviour of each part separately.
For large \( x \), Stirling's approximation for the Gamma function is given by:

$$
\Gamma(x) \sim \sqrt{2\pi x} \left( \frac{x}{e} \right)^{x}.
$$
We will apply this approximation to each Gamma term in both \( P(n) \) and \( F(n) \).
Each term in the product can be approximated using Stirling's formula. For large \( i \), we approximate:

$$
\Gamma(\beta i + 1) \sim \sqrt{2\pi (\beta i + 1)} \left( \frac{\beta i + 1}{e} \right)^{\beta i + 1},
$$
and
$$
\Gamma(\beta i + 1 - \alpha) \sim \sqrt{2\pi (\beta i + 1 - \alpha)} \left( \frac{\beta i + 1 - \alpha}{e} \right)^{\beta i + 1 - \alpha}.
$$
Thus, each ratio becomes

$$
\frac{\Gamma(\beta i + 1)}{\Gamma(\beta i + 1 - \alpha)} \sim \frac{\left( \beta i + 1 \right)^{\beta i + 1}}{\left( \beta i + 1 - \alpha \right)^{\beta i + 1 - \alpha}}\sim e^{-\alpha }(\beta i)^{\alpha } .
$$
Therefore, the product becomes

$$
P(n) \sim \prod_{i=1}^{n} e^{-\alpha }(\beta i)^{\alpha }  = e^{-\alpha n }(\beta^n\,n!)^{\alpha }\sim e^{-\alpha n }(\beta\,n)^{\alpha n }.
$$
Next, consider the final factor \( F(n) \). Using Stirling's approximation for large \( n \):

$$
\Gamma(\beta n + 1 - \alpha + \nu) \sim \sqrt{2\pi (\beta n + 1 - \alpha + \nu)} \left( \frac{\beta n + 1 - \alpha + \nu}{e} \right)^{\beta n + 1 - \alpha + \nu}.
$$
Thus, for large \( n \),

$$
F(n) \sim \frac{\left( \beta n + 1 - \alpha + \nu \right)^{\beta n + 1 - \alpha + \nu}}{\Gamma(1 - \alpha + \nu)} e^{-(\beta n + 1 - \alpha + \nu)}\sim e^{-\beta n }(\beta\,n)^{\beta n } .
$$
Combining the two asymptotic behaviours from \( P(n) \) and \( F(n) \), we get the overall behaviour of the generalised factorial. Ignoring constants and focusing on the dominant growth, we have the following {approximation},

$$
[n]_{\alpha,\beta,\nu}! \sim e^{-(\alpha+\beta) n }\cdot(\beta\,n)^{(\alpha+\beta) n }.
$$

\subsection{Carleman's Condition for Stieltjes Moment Problem}

{The Carleman condition provides a sufficient criterion for the uniqueness of the solution to the Stieltjes moment problem $\int_0^{\infty}x^n W(x)dx=m_n$ with $n=0,1,2...$.} It states that if the moments satisfy the following growth condition:
\begin{equation}
\sum_{n=1}^{\infty} m_n^{-\frac{1}{2n}} =\left\{\begin{array}{ll}
         \infty  \,\,\,\,\,\,\,\text{ the solution is unique}\\
        <\infty \,\,\,\text{non-unique solution may exist}.
        \end{array}\right.
\end{equation}
 Consider the moments defined {in the previous paragraph} by 

$$
m_n =e^{-(\alpha+\beta) n }\cdot(\beta\,n)^{(\alpha+\beta) n },
$$
where \(\alpha\) and \(\beta\) are constants. We will analyse Carleman's condition for the Stieltjes moment problem.

In our case, we can find:

$$
m_n^{-\frac{1}{2n}} = \left(e^{-(\alpha+\beta) n }\cdot(\beta\,n)^{(\alpha+\beta) n }\right)^{-\frac{1}{2n}} =e^{\frac{(\alpha+\beta) }{2} }\cdot(\beta\,n)^{-\frac{(\alpha+\beta) }{2} }.
$$
Now, we want to consider the series:

$$
\sum_{n=1}^{\infty} e^{\frac{(\alpha+\beta) }{2} }\cdot(\beta\,n)^{-\frac{(\alpha+\beta) }{2} }= e^{\frac{(\alpha+\beta) }{2}} \beta^{-\frac{(\alpha+\beta) }{2}}\cdot\sum_{n=1}^{\infty}  n^{-\frac{(\alpha+\beta) }{2} }.
$$
This series is a \( p \)-series and converges or diverges depending on the value of $\frac{(\alpha + \beta)}{2}$:

\begin{itemize}
    \item \textbf{Divergence}: If \(\frac{\alpha + \beta}{2} \leq 1\), the series diverges.
    \item \textbf{Convergence}: If \(\frac{\alpha + \beta}{2} > 1\), the series converges.
\end{itemize}

Finally, using Carleman's condition for the Stieltjes moment problem, we can summarise that:

\begin{itemize}
    \item If \(\frac{\alpha + \beta}{2} \leq 1\), then

    $$
    \sum_{n=1}^{\infty} m_n^{-\frac{1}{2n}} = \infty,
    $$
implies that the moment problem is determinate, and a unique positive measure exists and corresponds to the moments \( m_n \).
    
    \item If \(\frac{\alpha + \beta}{2} > 1\) the moment problem is indeterminate, a unique measure may not exist corresponding to the moments.
\end{itemize}

%--------------------------------

%------------------------------------------------------------------------------------------------------------------------------

\begin{thebibliography}{99}
%1
\bibitem{S26} Schr\"{o}dinger, E. (1926). Der stetige Übergang von der Mikro-zur Makromechanik. Naturwissenschaften, 14(28), 664-666.
%2
\bibitem{Glaub63} Glauber, R. J. (1963). The quantum theory of optical coherence. Physical Review, 130(6), 2529–2539. https://doi.org/10.1103/PhysRev.130.2529
%3
\bibitem{Klauder63} Klauder, J. R. (1963). Continuous‐representation theory. I. Postulates of continuous‐representation theory. Journal of Mathematical Physics, 4(8), 1055-1058.
%4
\bibitem{KPS2001}Klauder, J. R., Penson, K. A., Sixdeniers, J. M. (2001). Constructing coherent states through solutions of Stieltjes and Hausdorff moment problems. Physical Review A, 64(1), 013817.
%5
\bibitem{GBG2010} Gazeau, J. P., Baldiotti, M. C., Gitman, D. M. (2010). Coherent state quantisation and moment problem. Acta Polytechnica 50(3)
%6
\bibitem{KS91} Katriel, J., Solomon, A. I. (1991). Generalised q-bosons and their squeezed states. Journal of Physics A: Mathematical and General, 24(9), 2093.
%7
\bibitem{Herrmann2011} Herrmann, R. (2011). Fractional Calculus: An Introduction for Physicists. World Scientific
%8
\bibitem{KS85} Klauder, J.R., Skagerstam Be Sture, Coherent States, Application in Physics and Mathematical Physics (World Scientific, Singapore, 1985)
%8
\bibitem{PS99} Penson, K. A., Solomon, A. I., J. Math. Phys. 40, 2354 (1999)
%9
\bibitem{SP20002001}Sixdeniers, J-M., Penson, K. A., J.Phys. A. 33, 2907 (2000); 34, 2859 (2001)
%10
\bibitem{SPS99} Sixdeniers, J. M., Penson, K. A., Solomon, A. I. (1999). Mittag-Leffler coherent states. Journal of Physics A: Mathematical and General, 32(43), 7543.
%11
\bibitem{PBDHS2010} Penson, K. A., Blasiak, P., Duchamp, G.H.E., Horzela, A., Solomon, A. I. (2009). On certain non-unique solutions of the Stieltjes moment problem. \textit{Discrete Mathematics \& Theoretical Computer Science}, \textit{12}. 
%30
\bibitem{BDK94} Bonatsos, D., Daskaloyannis, C., Kolokotronis, P. (1994). Deformed oscillator algebras for two- dimensional quantum superintegrable systems. Phys. Rev. A 50, 3700.
%31
\bibitem{DGI2001} Draayer, J. P., Georgieva, A. I., Ivanov, M. I. (2001). Deformation of the boson sp(4, R) representation and its subalgebras. J. Phys. A: Math. Gen. 34, 2999.
%32
\bibitem{GGR95} Gueorguiev, V. G., Georgieva, A. I., Raychev, P. P. et al (1995). q-analog of $A_{m-1}\oplus A_{m-1}\subset A_{m\,n-1}$. Int. J. Theor. Phys. 34, 2095-2204. https://doi.org/10.1007/BF00673835.
%24
\bibitem{TP2007} Tavassoly, M. K., Parsaiean, A. (2007) Quantum statistical properties of some new classes of intelligent states associated with special quantum systems. J. Phys. A: Math. Theor. 40, 9905.
%25
\bibitem{RT2004} Roknizadeh, R., Tavassoly, M. K. (2004). The construction of some important classes of generalised coherent states: the nonlinear coherent states method. J. Phys. A: Math. Gen. 37, 8111.
%28
\bibitem{MMSZ97} Man'ko, V. I., Marmo, G., Sudarshan, E. C. G., Zaccaria, F. (1997). f-oscillators and nonlinear coherent states. Phys. Scr.55, 528.
%27
\bibitem{SIJ2009} Sadiq, M., Inomata, A., Junker, G. (2009). Coherent states for a polynomial SU(1,1) algebra and a conditionally solvable system. J. Phys. A: Math. Theor 42, 365210.
%29
\bibitem{BDK93} Bonatsos, D., Daskaloyannis, C., Kolokotronis, P. (1993). Generalised deformed SU(2) algebra. J. Phys. A: Math. Gen. 26, L871.
%12
\bibitem{GGM19} Garra, R., Giraldi, F.,  Mainardi, F. (2019). Wright-type generalised coherent states. WSEAS Trans. Math, 18, 428-431.
%13
\bibitem{GM23} Giraldi, F., Mainardi, F. Truncated generalised coherent states, J. Math. Phys. 64, 032105 (2023)
%20
\bibitem{Solomon94} Solomon, A. I. (1994). A characteristic functional for deformed photon phenomenology. Physics Letters A 196(1-2), 29-34.
%22
\bibitem{Dask91} Daskaloyannis, C. (1991). Generalised deformed oscillator and nonlinear algebras. J. Phys. A: Math. Gen. 24, L789.
%23
\bibitem{Dask92} Daskaloyannis, C. (1992). Generalised deformed oscillator corresponding to the modified Poschl-Teller energy spectrum. J. Phys. A: Math. Gen. 25, 2261.
%26
\bibitem{DY92} Daskaloyannis, C., Ypsilantis, K. (1992). A Deformed oscillator with Coulomb energy spectrum. J. Phys. A: Math. Gen. 25, 4157.
%14
\bibitem{Droghei21} Droghei, R. (2021). On a Solution of a Fractional Hyper-Bessel Differential Equation using a Multi-Index Special Function. Fract Calc Appl Anal 24, 1559–1570. https://doi.org/10.1515/fca-2021-0065
%15
\bibitem{A1965} Akhiezer, N. I. (1965). The classical moment problem and some related questions in analysis. Society for Industrial and Applied Mathematics.
%16
\bibitem{ST1943} Shohat, J., Tamarkin, J. (1943). The Problem of Moments (APS, New York)
%17
\bibitem{MSH09} Mathai, A. M., Saxena, R. K.,  Haubold, H. J. (2009). The H-function: theory and applications. Springer Science and Business Media.
%18
\bibitem{K2004} Kilbas, A. A. (2004). H-transforms: Theory and Applications. CRC Press.
%19
\bibitem{KS1994} Katriel, J., Solomon, A. I. (1994). Nonideal lasers, nonclassical light, and deformed photon states. Phys. Rev. A 49, 5149.
%21
\bibitem{MW95} Mandel, L., Wolf, E. (1995). Optical Coherence and Quantum Optics (Cambridge University Press, Cambridge).

\end{thebibliography}
\end{document}